\documentstyle[12pt]{article}
\renewcommand {\c}  {\'{c}}
\newcommand {\cc} {\v{c}}
\newcommand   {\s}  {\v{s}}
\begin{document}
\pagestyle{empty}
\vspace* {13mm}
\baselineskip = 24pt
\begin{center}

%
  {\bf   A NEW DEFORMED SUPERSYMMETRIC OSCILLATOR}
 \\[8.0mm]
S.Meljanac$^1$,  M.Milekovi\c$^{2,+}$ and A.Perica$^{1,++}$\\[7mm]
{\it $^1$ Rudjer Bo\s kovi\c \ Institute, Bijeni\cc ka c.54, 41001 Zagreb,
Croatia\\[3mm]
 $^2$ Prirodoslovno-Matemati\cc ki Fakultet, Zavod za teorijsku fiziku,
\\Bijeni\cc ka c.32, 41000 Zagreb, Croatia\\[3mm]
$^+$ e-mail: marijan@phy.hr\\
$^{++}$ e-mail: perica@thphys.irb.hr}
\bigskip

{\bf Abstract}

\end{center}

\vskip 0.2cm
%
We construct and discuss the Fock-space representation
for a deformed oscillator with "peculiar" statistics.
We show that corresponding algebra represents
deformed supersymmetric oscillator.
\begin{center}

{\it (to appear in Europhys.Lett.)}

\end{center}

%
%
%
%
\newpage
\setcounter{page}{1}
\pagestyle{plain}
\def\leer{\vspace{5mm}}
\baselineskip=24pt
\setcounter{equation}{0}%

{\it Introduction}. - The subject of quantum statistics,which is different from
the ordinary Bose and
Fermi statistics, has  attracted much attention in the past few years. One
motivation
comes from the study of some phenomena in condensed matter where the
dynamics is essentialy two-dimensional, thus allowing  anyon-like
statistics [1]. The other motivation comes from the theoretical and
experimental
search for possible violation of the Pauli exclusion principle in four
dimensions [2],
where quon-like statistics  [3] might play a significant role.\\
In either case, quantum groups and algebras [4] have offered a new
insight into the subject. The introduction of q-deformations of the
Heisenberg-Weyl
algebras has led to the investigation of particles interpolating between bosons
and fermions,  [5]. The q-bosons have been
introduced and discussed in a
variety of ways [6]. Particularly useful formulations of associative q-boson
algebras
are proposed through the Yang-Baxter R-matrix [7] which generalizes the notion
of permutational
 symmetry. The simplest such  algebras,associated to $4 \times 4$ R-matrices ,
were investigated to some extent in [8] and three types of deformed
algebras were found, among them a "peculiar" algebra which corresponded to the
R-matrix of the eight-vertex form.
 The detailed structure of  the "peculiar" algebra was not
discussed in [8] and it remains unclear whether this algebra makes a sense
physically, i.e. whether there exists a well-defined Fock-space representation
with positive norms and number operators.

In this Letter we investigate the structure of this "peculiar" algebra.
We construct and discuss the corresponding Fock-space representation
and show that norms of all states are positive definite.
We show that this algebra represents a new kind of deformed supersymmetric
oscillator.

\vskip 1cm

{\it "Peculiar" algebra}. - We start with the  following
( "peculiar" ) oscillator algebra [8]:
$$
(1-\epsilon q)a_{1}^{2} = \epsilon^{'}(1+\epsilon q)a_{2}^{2},\\[4mm]
$$
$$
a_{1}a_{2} = \epsilon^{''}a_{2}a_{1},\\[4mm]
$$
\begin{equation}
a_{1}a_{1}^{\dagger} = 1 +  (\epsilon q^{-1} +
 \frac{q^{-2}-1}{2})\; a_{1}^{\dagger}a_{1} + (\frac{q^{-2}-1}{2}) \;
a_{2}^{\dagger}a_{2},\\[4mm]
\end{equation}
$$
a_{2}a_{2}^{\dagger} = 1 +  ( - \epsilon q^{-1} +
\frac{q^{-2}-1}{2}) \; a_{2}^{\dagger}a_{2} + (\frac{q^{-2}-1}{2}) \;
a_{1}^{\dagger}a_{1},\\[4mm]
$$
$$
a_{1}a_{2}^{\dagger} = (\epsilon^{''}\frac{1+q^{-2}}{2}) \;
a_{2}^{\dagger}a_{1} +
(\epsilon^{'}\frac{q^{-2}-1}{2}) \; a_{1}^{\dagger}a_{2} ,\\
$$
$$
a_{2}a_{1}^{\dagger} = (\epsilon^{''}\frac{1+q^{-2}}{2}) \;
a_{1}^{\dagger}a_{2} +
(\epsilon^{'}\frac{q^{-2}-1}{2}) \; a_{2}^{\dagger}a_{1} ,
$$
where $q\in {\bf R}$, $\epsilon^{2}=\epsilon^{'2}=\epsilon^{''2}=1$.
When $q^{2}=1$, the above algebra, eq.(1), represents one Bose and one
Fermi oscillator which commute or anticommute (depending on whether
$\epsilon^{''}$ is 1 or -1).
We observe that the "peculiar" algebra, eq.(1), has no well-defined
number operators $N_{1}$, $N_{2}$, in the usual sense:
$[N_{i},a_{j}]=-a_{i}\delta_{ij}$, $[N_{i},a_{j}^{\dagger}]=a_{i}\delta_{ij}$,
$i,j=1,2$.
{}From $[N_{1},a_{1}]=-a_{1}$, it follows that $[N_{1},a_{1}^{2}]=-2a_{1}^{2}$.
Owing to eq.(1) one obtains $[N_{1},a_{2}^{2}]=-2a_{2}^{2}$, which
contradicts the demanded relation $[N_{1},a_{2}]=0$.
Hence the relations $[N_{1},a_{2}]=[N_{2},a_{1}]=0$ contradict eq.(1).
However,the total number operator $N=N_1 +N_2$ is well defined.Of course, when
$q^{2}=1$,
 the number operators $N_1$ and $N_2$ are also well defined, i.e.
$N_{1,2}=N_{B,F}$ .

{\it Fock-space representation of "peculiar" algebra}. - Let us assume that
there is a vacuum $|0,0>$ satisfying
$a_{i}|0,0> = 0$,  $i=1,2$.
 The excited states can be constructed by multiple action of the
 operators $a_{1}^{\dagger}$ and $a_{2}^{\dagger}$ on the vacuum $|0,0>$ and
 are of the form
\begin{equation}
|n_{1},n_{2}> \; \propto
\;(a_{1}^{\dagger})^{n_{1}}(a_{2}^{\dagger})^{n_{2}}|0,0>,
\qquad n_{1},n_{2}\in {\bf N}.
\end{equation}
Notice that the action of $N_1\, (N_2)$ on the states (2) is not well defined
for
$n_2\geq 2 \\ (n_1 \geq 2)$ and hence, in general, $n_1\, (n_2)$ are not
eigenvalues of
$N_1\, (N_2)$.This is the consequence of the quadratic relation $a_1^2 \propto
a_2^2$
 (eq.(1)) for $q^2 \neq 1$.
Furthermore, we find that the $|n_{1},n_{2}>$ states are degenerate
(linearly dependent) for the fixed sum $n_{1}+n_{2}=n$, in the  following
sense:
$|n_1,n_2> \propto |n_1-2k,n_2+2k>$,$\;$ $k \in {\bf Z}$,\\$n_1-2k \geq
0$,$n_2+2k\geq 0$ and
$|n_1+1,n_2-1> \propto |n_1+1-2k,n_2-1+2k>$,\\$k \in {\bf Z}$,$n_1+1-2k \geq
0$,
$n_2-1+2k\geq 0$.The states for fixed n
can be reduced to two states , $|n,0>$ and
$|(n-1),1>$ or alternatively ,to $|0,n>$ and $|1,(n-1)>$.\\
Hence,the complete set of states can be represented by two symmetric pictures
(for $q^2\neq 1$)
\begin{equation}
|n,\nu>\; \propto \; (a_{1}^{\dagger})^{n}(a_{2}^{\dagger})^{\nu}|0,0>,\;\qquad
(a)
\end{equation}
$$
|\nu,n>\; \propto \; (a_{1}^{\dagger})^{\nu}(a_{2}^{\dagger})^{n}|0,0>,\;\qquad
(b)
$$
where $n \in {\bf N_0}$, $\nu = 0,1$.Now, n $\,(\nu)$ is eigenvalue of $N_1 \,
(N_2)$
in the picture (3a) i.e.  $N_2 \, (N_1)$ in the picture (3b).In the following
we use the
first picture (3a).\\  There are two towers of states
generated by the $a_{1}^{\dagger}$ creation operator. One tower is $|n,0>$,
generated from the $|0,0>$ vacuum, ($\nu = 0$). The other tower is $|n,1>$,
generated from the second vacuum $|0,1>$, ($\nu = 1$).
Using the algebra (1) we find that
$$
a_{1}^{\dagger}a_{1}(a_{1}^{\dagger})^{n}(a_{2}^{\dagger})^{\nu}|0,0> =
\phi_{1}(n,\nu)(a_{1}^{\dagger})^{n}(a_{2}^{\dagger})^{\nu}|0,0>,\\
$$
\begin{equation}
\phi_{1}(n,\nu) = \frac{1}{2}[n]_{(\epsilon q)^{-1}}(1 + (\epsilon q)^{
1-n-2\nu}),
\end{equation}
where $[n]_{(\epsilon q)^{-1}} = \frac {(\epsilon q)^{-n} - 1}
{(\epsilon q)^{-1} - 1}$, and $n \in {\bf N_0}$, $\nu = 0,1$.\\
It is important to observe that, for $q \in {\bf R}$, the function $\phi_1$
is positive :
$\phi_{1}(n,\nu) > 0$, $\forall n\in {\bf N_0}$, $\nu = 0,1$.
Furthermore, $\phi(n,\nu)$ cannot be written as a function of one variable.
If it could be done so, this would mean that there would be only one tower
of states, and that $a_1 \propto a_2$.
Hence, all states $(a_{1}^{\dagger})^{n}(a_{2}^{\dagger})^{\nu}|0,0>$, eq.(3a),
have positive definite norms and can be normalized. The normalized states are\\
\begin{equation}
|n,\nu> = \frac{(a_{1}^{\dagger})^{n}(a_{2}^{\dagger})^{\nu}|0,0>}
{\sqrt{[\phi_{1}(n,\nu)]!}},\\
\end{equation}
where $[\phi_{1}(n,\nu)]! = \phi_{1}(n,\nu)\cdots \phi_{1}(1,\nu)$,
$[\phi_{1}(0,\nu)]! = 1$, and the orthonormality condition reads
$<n,\nu|n^{'},\nu^{'}> = \delta_{nn^{'}}\delta_{\nu \nu^{'}}$,
$ \nu ,\nu^{'}=0,1$.
Owing to this orthonormality relation, any linear combination
of states, eq.(5), has a positive norm.
Specially,
$$
||\alpha |n,0> + \beta |n-1,1>||^{2} = |\alpha |^{2} + |\beta |^{2} > 0.
$$
It is easy to find the action of the $a_{i}$,$a_{i}^{\dagger}$ operators
on the states ,eq.(5),namely:
\bigskip
\begin{equation}
a_{1}^{\dagger}|n,\nu> = \sqrt{\phi_{1}(n+1,\nu)}|n+1,\nu>,\\[4mm]
\end{equation}
$$
a_{1}|n,\nu> = \sqrt{\phi_{1}(n,\nu)}|n-1,\nu>,\\[4mm]
$$
$$
a_{2}^{\dagger}|n,\nu> = \sqrt{\frac{[\phi_{1}((n+2\nu),(1-\nu))]!}
{[\phi_{1}(n,\nu)]!}}(\frac{1 - \epsilon q}{\epsilon^{'}(1 + \epsilon q)}
)^{\nu}(\epsilon^{''})^{n}|(n+2\nu),(1-\nu)>,\\ [4mm]
$$
$$
a_{2}|n,\nu> = \sqrt{\frac{[\phi_{1}(n,\nu)]!}{[\phi_{1}((n-2+2\nu),(1-\nu)
)]!}}(\frac{1 - \epsilon q}{\epsilon^{'}(1 + \epsilon q)})^{1-\nu}(\epsilon^
{''})^{n}|(n-2+2\nu),(1-\nu)>.
$$
In the picture (3a), the $a_{1}^{\dagger}$ operator builds two infinite towers
on $|0,0>$ and $|0,1>$, respectively, whereas the $a_{2}$, $a_{2}^{\dagger}
$ operators interconnect the two towers.\\
In the  picture (3b), in which the indices are interchanged,$1\leftrightarrow
2$ and \\$\varepsilon \leftrightarrow - \varepsilon $, the $a_{2}^{\dagger}$
operator creates two towers based on $|0,0>$ and $|0,1>$ while the
$a_{1}^{\dagger}$ operator braids between these two towers.All equations (5-6)
hold with $1\leftrightarrow
2$ , $\varepsilon \leftrightarrow - \varepsilon $.

{\it Deformed SUSY oscillator}. - We can define the operators $Q_{ij} =
a_{i}a_{j}^{\dagger}$, ($Q_{ij}^
{\dagger} = Q_{ji}$) and $\tilde{Q_{ij}} = a_{i}^{\dagger}a_{j}$,
($\tilde{Q_{ij}}^{\dagger} = \tilde{Q_{ji}}$), $i,j = 1,2$, satisfying
(in the picture (3a))
\begin{equation}
\begin{array}{c}
Q_{ij}= \delta_{ij} + p'\; R_{ki,jl} \, \tilde{Q}_{kl},\\[4mm]
[N,Q_{ij}] = [N,\tilde{Q_{ij}}] = 0,\qquad  \forall i,j=1,2,\\[4mm]
Q_{11}|n,\nu> = \phi_{1}(n+1,\nu)|n,\nu>,\\[4mm]
Q_{22}|n,\nu> = \phi_{2}(n+2\nu,1-\nu)|n,\nu>,\\[4mm]
Q_{12}|n,\nu> = \psi_{12}(n,\nu)|n-1+2\nu, 1-\nu>,\\[4mm]
Q_{21}|n,\nu> = \psi_{21}(n,\nu)|n-1+2\nu, 1-\nu>,\\[4mm]
Q_{12}^{\dagger}Q_{12} = Q_{21}Q_{12} = \psi_{12}^{2}(n,\nu),\\[4mm]
Q_{12}^{2} = \psi_{12}(n,\nu)\psi_{12}(n-1+2\nu,1-\nu) ,
\end{array}
\end{equation}
where
$$
\phi_{2}(n,\nu) = \frac{[\phi_{1}(n,\nu)]!}{[\phi_{1}(n-2+2\nu,1-\nu)]!}
(\frac{1-\epsilon q}{\epsilon ' (1+\epsilon q)})^{2(1-\nu)},\\[4mm]
$$
\begin{equation}
\psi_{12}(n,\nu)=\sqrt{\frac{\phi_{1}(n-2+2\nu,1-\nu)
[\phi_{1}(n+2\nu,1-\nu)]! }{[\phi_{1}(n,\nu)]!}} (\frac{1-\epsilon q}
{\epsilon ' (1+\epsilon q)})^{\nu} (\epsilon ^{''} )^{n}\\[4mm]
\end{equation}
$$
\psi_{21}(n,\nu) = \psi_{12}(n-1+2\nu,1-\nu).\\
$$
Analogous relations can be obtained for the operators $\tilde{Q}_{ij}$
 and  $\tilde{Q}_{ij}^{\dagger}$ using eqs.(6,7).
Notice that $(Q_{12})^{2} \neq 0$ $((\tilde{Q}_{12})^2 \neq 0)$ when $q^{2}
\neq 1$, and
 $(Q_{12})^{2} = 0$ $((\tilde{Q}_{12})^2 = 0)$ when $q^{2} = 1$.\\
We can define the Hamiltonian H as
\begin{equation}
\begin{array}{l}
\{ Q_{12},Q_{12}^{\dagger}\} = 2H,\\[4mm]
[H,Q_{12}] = [H,Q_{12}^{\dagger} ] = [H,N] = 0,\\[4mm]
H | n,\nu> = \frac{1}{2} ( \;(\psi_{12}(n,\nu))^2 +
(\psi_{12}(n-1+2\nu,1-\nu))^2 \;)| n,\nu>,
\end{array}
\end{equation}
and similarly, the Hamiltonian $\tilde{H}$ as\\
\begin{equation}
\begin{array}{l}
\{ \tilde{Q}_{12},\tilde{Q}_{12}^{\dagger}\} = 2\tilde{H},\\[4mm]
[\tilde{H},\tilde{Q}_{12}] = [\tilde{H},\tilde{Q}_{12}^{\dagger} ] =
[\tilde{H},N] = 0,\\[4mm]
\tilde{H} | n,\nu> = \frac{1}{2} ( \;(\tilde{\psi}_{12}(n,\nu))^2 +
(\tilde{\psi}_{12}(n-1+2\nu,1-\nu))^2 \;)| n,\nu>.\\
\end{array}
\end{equation}
The relations (9) and (10) define a new kind of q-deformation of the
supersymmetric
(SUSY) oscillator [9].
We point out that the spectrum of H ($\tilde{H}$) is positive and degenerate,
i.e.
the states $|n,0>$ and $|n-1,1>$ have the same energy
$\frac{1}{2}((\psi(n,0))^2 + (\psi(n-1,1))^2 )$.These properties are  typical
for SUSY oscillator,except that for $q^2 \neq 1$ the energy levels are not
equidistant.In the limit $q=+1$,the state $|n,0>$ ($|n-1,1>$) is bosonic
(fermionic) in the picture (3a).In the limit $q=-1$,the state $|0,n>$
($|1,n-1>$) is bosonic (fermionic)in the picture (3b).\\
 The q-deformed SUSY algebra  (9) is generated by the set
$\{ N,Q_{12},Q_{12}^{\dagger},H \}$ and the q-deformed SUSY algebra  (10)  by
the set
$\{ N,\tilde{Q}_{12},\tilde{Q}_{12}^{\dagger},\tilde{H} \}$.Notice that our
Hamiltonian H (and $\tilde{H}$) is invariant under the q-superalgebra since
H and Q ($\tilde{H}$ and $\tilde{Q}$) mutually commute, in contrast to
the Hamiltonian of the form $H=\{Q_{+},Q_{-}\}$ mentioned in [10].
The q-deformed supercharges,operators $Q_{ij}$, $\tilde{Q}_{ij}$, $i \neq j$,
also braid between
the two towers and preserve the total number operator $N = N_{1} + \nu$.
Although the operators Q and $\tilde{Q}$ are not nilpotent ($Q^2_{12}\neq 0 $
for
$q^2 \neq 1$,contrary to the ordinary SUSY oscillator),their irreducible
representations remain two-dimensional,as a consequence of the relation $a^2_1
\propto a^2_2 $ (eq.(1)).

{\it Conclusion}. - In conclusion,we have constructed and investigated the
Fock-space
representation
  for the "peculiar" algebra defined for a two-mode
oscillator in ref.[8]. We have shown that this algebra corresponds to the
deformed
supersymmetric oscillator.This deformed SUSY oscillator represents an
alternative
mechanism for the violation of the Pauli exclusion principle [2].It is also
possible to generalize the
"peculiar" algebra to a multimode case and to include arbitrary
relations between powers of the operators $a_{i}$ with arbitrary exponents.

\newpage

{\bf REFERENCES}

\begin{description}
\item{[1]}
Leinaas J.M.and Myrheim J.,
{\it Nuovo Cim.}, {\bf 37} (1977) 1\\
Wilczek F.,
{\it Phys.Rev.Lett.}, {\bf 48} (1982) 1144;ibid., {\bf 49} (1982) 957;\\
Wu Y.S.,
{\it Phys.Rev.Lett.},{\bf 52} (1984) 111;ibid., {\bf 53} (1984) 111.
\item{[2]}
Ignatiev A.Yu. and Kuzmin V.A., {\it Sov.J.Nucl.Phys.}, {\bf 46} (1987) 444;\\
Mohapatra R.N.and Greenberg O.W., {\it Phys.Rev.D}, {\bf  39} (1989) 2032;\\
Miljani\c $\,$ D. $\;$ et al.,{\it Phys.Lett.B}.,{\bf  252} (1990) 487;\\
Ramberg E. and Snow G.A.,{\it Phys.Lett.B},{\bf  242} (1990) 407.
\item{[3]}
Greenberg O.W.,
{\it Phys.Rev.Lett.},{\bf 64} (1990) 705;
{\it Phys.Rev.D},{\bf 43} (1991) 4111;\\
 Meljanac S. and Perica A.,
{\it Mod.Phys.Lett. A},{\bf 9} (1994) 3293.
\item{[4]}
 Drinfeld V.G.,
in {\it Proceedings of the ICM}  (Berkeley, CA ,1986, p.798);\\
 Jimbo M.,
{\it Lett.Math.Phys.},{\bf 10} (1985) 63.
\item{[5]}
Biedenharn L.C.,
{\it J. Phys.A :Math.Gen.},{\bf 22} (1989) L873;\\
Macfarlane A.J.,
{\it J. Phys.A :Math.Gen.},{\bf 22} (1989)  L983.
\item{[6]}
Tuszynski J.A.,Rubin J.L.,Meyer J. and Kibler M.,
{\it Phys.Lett.A},{\bf 175}  (1993) 173;\\
Meljanac S.,Milekovi\c $\,$M.  and Pallua S.,
{\it Phys.Lett.B},{\bf 328}  (1994) 55;\\
 Bardek V.,  Dore\s i\c $\;$M. and  Meljanac S.,
 {\it  Int.J.Mod.Phys.A},{\bf 9} (1994) 4185 ;\\
Bardek V.,  Dore\s i\c $\;$M. and  Meljanac S.,
{\it Phys.Rev.D },{\bf 49} (1994) 3059;\\
Bonatsos D. and Daskaloyannis C.,
{\it Phys.Lett.B},{\bf 278} (1992) 1;\\
Odaka K.,Kishi K. and Kamefuchi S.,
{\it J. Phys.A :Math.Gen.},{\bf 24} (1991) L591.
\item{[7]}
 Kulish P.P.,
{\it Phys.Lett.A}, {\bf 161}  (1991) 50;\\
 Fairlie D. and C. Zachos C.,
{\it Phys.Lett.B}, {\bf 256} (1991) 43;\\
Meljanac S.,Milekovi\c $\,$ M. and Perica A.,
{\it Europhys.Lett.},{\bf 28} (1994) 79.
\item{[8]}
 Van der Jeugt J.,
{\it J. Phys.A :Math.Gen.}, {\bf 26} (1993) L405.
\item{[9]}
de Crombrugghe M. and  Rittenberg V.,
{\it Ann.Phys.}, {\bf 151} (1983) 99;\\
 Chaichian M. and  Kulish P.,
{\it Phys.Lett.B}, {\bf 234} (1990) 72.
\item{[10]}
 Parthasarathy R. and Wiswanathan K. S.,
{\it J. Phys.A :Math.Gen.}, {\bf 24} (1991) 613.

\end{description}

\end{document}